\documentclass[10pt]{amsart}
\usepackage{amscd}
\usepackage{epsfig}

\newenvironment{theorem}{\medskip \noindent\normalsize {\sc Theorem.}\it}{\medskip}
\newenvironment{lemma}{\medskip \noindent\normalsize {\sc Lemma.}\it}{\medskip}
\newenvironment{prop}{\medskip \noindent\normalsize {\sc Proposition.}\it}{\medskip}
\newenvironment{corollary}{\noindent\normalsize {\sc Corollary.}\it}{\medskip}

\def\H{{\mathcal H}}
\def\A{{\mathcal A}}
\def\D{{\mathcal D}}
\def\J{{\mathcal J}}
\def\L{{\mathcal L}}
\def\R{{\mathcal R}}

\def\Z{{\mathbb Z}}

\def\dim{{\rm dim}}

\begin{document}
\title{Cells and representations of right-angled Coxeter groups}
\author{Mikhail Belolipetsky}
\address{Mikhail Belolipetsky}
\address{Max Planck Institute of Mathematics, Vivatsgasse 7, 53111 Bonn, Germany}
\address{Sobolev Institute of Mathematics, Koptyuga 4, 630090 Novosibirsk, Russia}
\email{mbel@math.nsc.ru}
\date{}

\begin{abstract}
We study Kazhdan--Lusztig cells and the corresponding representations of right-angled Coxeter
groups and Hecke algebras associated to them. In case of the infinite groups generated
by reflections in the hyperbolic plane about the sides of right-angled polygons we obtain an
explicit description of the left and two-sided cells. In particular, we prove that
there are infinitely many left cells but they all form only three two-sided cells.
\end{abstract}

\maketitle

\section{Introduction}

A Coxeter group is said to be right-angled if for any two distinct simple reflections
their product has order $2$ or $\infty$. Due to their special properties and rich
structure, right-angled Coxeter groups often arise in different geometric and algebraic
problems. In a sense, the most interesting right-angled Coxeter groups are those which
can be presented as groups generated by reflections in hyperbolic spaces. Examples of
geometric applications of such groups can be found in \cite{Davis}, \cite{MV}, \cite{BGI}.
These groups also occur as Weyl groups of certain Kac--Moody
Lie algebras \cite{KP}. We are interested in the representations of right-angled
Coxeter groups $W$ and the corresponding Hecke algebras $\H$.

An important difference between hyperbolic reflection groups and affine or
finite Coxeter groups is that the former do not usually have a local system of generators
in the sense of \cite{OV}. In particular, the groups $P_n$ (see~\ref{Pn}) that are
generated by reflections in the hyperbolic plane about the sides of right-angled \mbox{$n$-gons}
can be called anti-local since only adjacent generators in the canonical system of
simple reflections commute. This means that we can hardly hope to construct the
representations of these groups inductively using the approach suggested in~\cite{OV}.
At the same time, we can still use the global methods of \cite{KL} and try to describe the
Kazhdan--Lusztig cells in our groups. It appears that certain symmetries of the initial
groups are reflected in the structure of the partitions into cells making the
cells trackable. Our main results concern the groups $P_n$ but the methods can be
applied to the other right-angled Coxeter groups as well.

We show that the partition of the group $P_n$ ($n \ge 5$) into left cells consists of
infinitely many elements and give an explicit description of the cells. At the same time,
all the left cells in our groups form precisely three two-sided cells (one of which is,
of course, the trivial cell). It was first
shown in \cite{Bedard} that the number of one-sided cells of a hyperbolic Coxeter group
can be infinite, but there the author used an implicit argument and obtained
only a conjectural structure of the corresponding cell-partitions (see also \cite{Cass}
for discussion). The fact that infinitely many left cells can still fall into finitely
many two-sided equivalence classes and hence give rise to only finitely many $\H$-bimodules
seems to have been previously unnoticed.

The explicit description of the cells makes it possible to consider corresponding
representations of the Coxeter groups and their Hecke algebras. Here we only start
the related analysis leaving more detailed considerations for the future. We discuss
the representations using $W$-graphs which were also introduced in \cite{KL}.

\smallskip

While working on this paper, I enjoyed the hospitality of the MPIM in Bonn. I~wish to
thank A.~Vershik for several helpful discussions. V.~Ostrik read an early
version of the paper and gave me some important suggestions. Finally I am grateful to
W.~Casselman for his response on my submission to arXiv and lots of email conversations.

\section{Preliminaries}
We recall some well-known facts about Coxeter groups and Hecke algebras associated
with them. The basic reference is the fundamental paper \cite{KL}. All the material
cited here can be also found in the book \cite{H}.

\subsection{} Let $W$ be a Coxeter group and let $S$ be the corresponding set of
simple reflections. With some ambiguity of language we shall also call by Coxeter
group the Coxeter system $(W,S)$. The {\it Hecke algebra} $\H$ over the ring
$\A = \Z[q^{1/2}, q^{-1/2}]$ of Laurent polynomials in $q^{1/2}$ is defined as
follows. As an $\A$-module, $\H$ is free with basis $T_w$ ($w\in W$), the
multiplication is defined by
\begin{align*}
& T_wT_{w'} = T_{ww'},\ {\rm if}\ l(ww') = l(w)+l(w'),\\
& T_s^2 = q + (q-1)T_s,\ {\rm if}\ s\in S,
\end{align*}
where $l(w)$ is the length of $w$ in $(W,S)$.

It will be also convenient to define
$$ \widetilde{T}_w = q^{-l(w)/2}T_w.$$

\subsection{} Let $a\to \overline a$ be the involution of the ring $\A$ defined by
$\overline{q^{1/2}} = q^{-1/2}$. This extends to an involution $h\to\overline h$ of the
ring $\H$ given by
$$\overline{\sum a_wT_w} = \sum \overline{a}_w T_w^{-1}.$$

Let $\le$ be the Bruhat order on $W$ (see \cite{H}, Ch.~5.9). Denote, as usual,
$q_w = q^{l(w)}$, $\epsilon_w = (-1)^{l(w)}$ for all $w\in W$.
In \cite{KL} it was shown that for any $w\in W$ there exists a unique element
element $C_w\in \H$ such that
\begin{align*}
& \overline {C}_w = C_w,\\
& C_w = \epsilon_wq_w^{1/2}\sum\limits_{y\le w}\epsilon_y q_y^{-1} \overline{P_{y,w}}T_y,
\end{align*}
where $P_{y,w}\in\A$ is a polynomial in $q$ of degree $\le\frac12(l(w)-l(y)-1)$ for
$y < w$, and $P_{w,w} = 1$.

Elements $C_w$ form a basis (called {\it $C$-basis}) of $\H$ as an $\A$-module. This basis
and the polynomials $P_{y,w}$ (called {\it Kazhdan--Lusztig polynomials}) turn out to be
of fundamental interest in the representation theory of Coxeter groups and Hecke
algebras.

\subsection{} \label{2.3}
Given $y,w\in W$, we write $y\prec w$ if $y < w$ and $P_{y,w}$ is a polynomial
in $q$ of degree exactly $\frac12(l(w)-l(y)-1)$ (which is, of course, possible only if
$l(w)-l(y)$ is odd). In this case the coefficient of the highest power of $q$ in $P_{y,w}$
is denoted by $\mu(y,w)$. If $w\prec y$ we set $\mu(w,y) = \mu(y,w)$, otherwise (if neither
$y \prec w$ nor $w \prec y$) let $\mu(y,w) = \mu(w,y) = 0$. We write $y - w$ if
$\mu(x,y) \neq 0$.

For any $w\in W$ define subsets of $S$:
$$ \L(w) = \{s\in S \mid sw < w \},\quad \R(w) = \{s\in S \mid ws < w \}. $$

Now define $y\le_L w$ to mean that there is a chain $y = y_0, y_1,\dots,y_n=w$ such that
$y_i-y_{i+1}$ and $\L(y_i)\not\subset \L(y_{i+1})$ for $0\le i < n$. Similarly, say that
$y\le_R w$ if $y^{-1} \le_L w^{-1}$. Finally, define $y\le_{LR}w$
to mean that there exists a chain $y = y_0, y_1,\dots,y_n=w$ such that for each $i < n$
either $y_i\le_L y_{i+1}$ or $y_i\le_R y_{i+1}$.

Let $\sim_L$, $\sim_R$, $\sim_{LR}$ be the equivalence relations associated to the
preorders $\le_L$, $\le_R$, $\le_{LR}$, respectively. The corresponding equivalence
classes are called {\it left}, {\it right} and {\it two-sided cells} of $W$.

\subsection{} We shall often make use of the following properties of the defined
relations.

\label{lemma1}\begin{lemma}
Let $x,y\in W$ and $x < y$.
\begin{itemize}
\item[(i)] If there exists $s\in S$ such that $x < sx$, $sy < y$, then $x\prec y$
if and only if $y = sx$.
\item[(ii)] If there exists $s\in S$ such that $x < xs$, $ys < y$, then $x\prec y$
if and only if $y = xs$.
\item[] Moreover, in each of the cases $\mu(x,y) = 1$.
\end{itemize}
\end{lemma}

\noindent (This statement can be found in the proof of Theorem~1.3 in \cite{KL};
it follows from the formula for the action of the elements $T_s$ on the
basis $\{C_w \mid w\in W\}$ of $\H$ and the definition of $C_w$.)

\begin{corollary}\

\begin{itemize}
\item[(i)] If $x \le_L y$, then $\R(x)\supset \R(y)$. Hence, $x\sim_L y$ implies $\R(x) = \R(y)$.
\item[(ii)] If $x \le_R y$, then $\L(x)\supset \L(y)$. Hence, $x\sim_R y$ implies $\L(x) = \L(y)$.
\end{itemize}
\end{corollary}

\noindent (To prove the corollary it is enough to consider the case $x - y$ with
$\L(x)\not\subset \L(y)$, details can be found in \cite{KL}.)

\subsection{} The main purpose of the partition of a Coxeter group $W$ into cells
is that it gives rise to the representations of group $W$ and its Hecke algebra $\H$.
It is convenient to describe these representations using $W$-graphs.

Let $X$ be a set. Consider an oriented graph $\Gamma$ whose set of vertices is
$X$; for each $x\in X$ there is assigned a subset $I_x$ of $S,$ and if
$I_x \not\subset I_y$, then there is an edge $(x,y)\in X\times X$ labeled by
an integer $\mu(x,y).$ Graph $\Gamma$ is called {\it a $W$-graph} if the map
$s\to\tau_s$, such that
$$
\tau_sx = \left\{ \begin{array}{l}
 -x,\quad {\rm if}\ x\in X,\ s\in I_x, \phantom{\displaystyle\sum\limits_x} \\
qx + q^{1/2}\displaystyle\sum\limits_{y\in X,\ s\in {I_y}}\mu(x,y)y,\quad
{\rm if}\ x\in X,\ s\not\in I_x,
\end{array}\right.
$$
defines a representation of $\H$ on the free $\A$-module $\A(X)$.

In \cite{KL} it was shown that $X = W$ with $I_x = \L(x)$ and $\mu(x,y)$ defined by
the polynomial $P_{x,y}$ as in \ref{2.3} gives a $W$-graph.
Moreover, its full subgraphs corresponding to the left cells with the same sets $I_x$
and the same function $\mu$ are $W$-graphs themselves.
Finally, when $W$ is a symmetric group $S_n$ it was proved that all the irreducible
representations of $\H$ are defined by the $W$-graphs associated to the left cells of $W$.

\section{Right-angled Coxeter groups}
\subsection{} \label{Pn}
Recall that Coxeter group $(W,S)$ is called {\it a right-angled Coxeter group}
if for any $s\neq t$ in $S$ the product $st$ has order $2$ or $\infty$.
The Coxeter graph of $W$ has only edges labeled via $\infty$ (see examples in Figure~1).
\begin{figure}[ht]
\vskip .5cm
\psfig{file=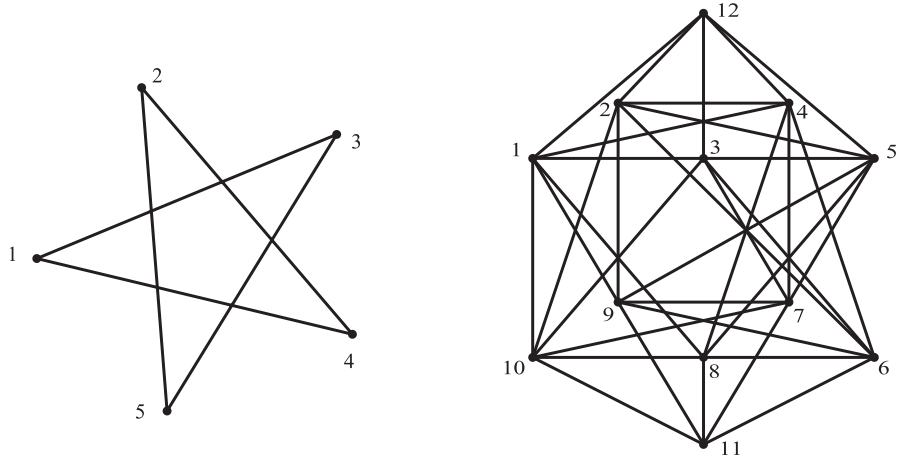, scale=.80}
\caption{Coxeter graphs of the groups generated by reflections about the sides
of a right-angled pentagon in the hyperbolic plane (the group $P_5$) and the
regular right-angled dodecahedron in the hyperbolic $3$-space.}
\end{figure}

The simplest example of a right-angled Coxeter group is the infinite dihedral group
$D_\infty =\ <s_1, s_2 \mid s_1^2 = s_2^2 = 1>$. We are mainly interested in right-angled
Coxeter groups generated by reflections in the hyperbolic $n$-space.
The important example for us are the groups $P_n$ having representations
$$P_n =\ <s_1, s_2,\dots,s_n \mid s_i^2 = 1,\ (s_js_{j+1})^2 = 1,\ (s_ns_1)^2 = 1>$$
where $i = 1,\dots ,n;$ $j = 1, \dots n-1$. If $n\geq 5$ the group $P_n$ can be
presented as a group of isometries of the hyperbolic plane generated by reflections
about the sides of a right-angled hyperbolic $n$-gon (for more information about this and
other geometric facts mentioned here we refer the reader to \cite{Vinb}).

\subsection{} \label{sec32}
If a Coxeter group $(W,S)$ is presented by isometries of the hyperbolic space, then it
gives rise to a tessellation of the space by the fundamental chambers of the group.
The dual graph $G$ of this tessellation endowed with the standard graph metric reflects
the structure of the initial group: fixing a vertex $e\in G$ which will correspond
to the identity element of $W$ and labelling all the edges with the corresponding
simple reflections from $S$, we can associate to an element $w\in W$ represented by a
word on $S$ a geodesic path starting from $e$ in $G$. Two paths define the same element if
and only if they have the same ends. Reduced expressions in $W$ correspond to the shortest
geodesics in $G$. The graph $G$ is called a Cayley graph of $(W,S)$; this graph can
be defined for an arbitrary group with a fixed system of generators.

\subsection{} \label{linebasics}
Let $(W,S)$ be an arbitrary Coxeter group. We call by {\it lines} elements $w\in W$ that
have unique reduced expressions. By a subword of a word $s_1s_2\dots s_n$, $s_i\in S$
we mean any expression of the form $s_is_{i+1}\dots s_j$, $1\le i \le j \le n$. The subword
$u$ of word $w_0uw_1$ is called {\it a segment} of the word and of the corresponding
element $w\in W$ if $u$ represents a maximal line such that any reduced expression
of $w$ has the form  $w_0'uw_1'$ with
$w_0 = w_0'$, $w_1 = w_1'$ as the elements of the group $W$ (here ``maximal'' means
that $u$ is not contained in any other subword with the same properties).

Lines are exactly the unique shortest geodesics of a Cayley graph that start from $e$.
The same way, segments of $w\in W$ correspond to the (maximal) geodesics of a Cayley graph
that are contained in any shortest geodesic corresponding to $w$. The only segment of a line
is the line itself.

Let us give a characterization of the lines and segments in a right-angled Coxeter
group.

\begin{prop} Let $(W,S)$ be a right-angled Coxeter group. Then
\begin{itemize}
\item[(i)] a word $s_1\dots s_k$ represents a line in $W$ if and only if for any $i=1,\dots ,k-1$
the product $s_is_{i+1}$ has order $\infty$;
\item[(ii)] a subword $u$ of $w_0uw_1$ is a segment if and only if $u$ gives a minimal
line such that $\#\R(w_0)\neq 1$ and  $\#\L(w_1)\neq 1$ (here ``minimal'' means that
$u$ does not contain any proper or empty subword with the same property).
\end{itemize}
\end{prop}

The proof easily follows from the definitions.
\medskip

The proposition implies that any element $w$ in a right-angled Coxeter group can be written
as $w = u_0s_{1,1}\dots s_{1,n_1}u_1 s_{2,1}\dots s_{2,n_2}u_2 \dots  s_{k,1}\dots s_{k,n_k}u_k$
where all the subwords $u_i$ are either trivial or segments in $w$, $n_i\ge 2$ and
$s_{i,j}s_{i,j+1} = s_{i,j+1}s_{i,j}$ for $i = 1,\dots,k$, $j = 1,\dots, n_i-1$.

\section{Distinguished Involutions}
\vbox{
We recall some definitions and results from \cite{L1a}, \cite{L2a}.

As in Section~2 we fix a Coxeter group $(W,S)$ and denote by $\H$ the corresponding
Hecke algebra over $\A = \Z[q^{1/2}, q^{-1/2}]$. For any $x,y,z\in W$ define elements
$f_{x,y,z}$, $g_{x,y,z}$, $h_{x,y,z}$ in $\A$ so that
}
\begin{align*}
\widetilde{T}_x\widetilde{T}_y &= \sum\limits_{z}f_{x,y,z}\widetilde{T}_z,\\
\widetilde{T}_xC_y &= \sum\limits_{z}g_{x,y,z}C_z,\\
C_xC_y &= \sum\limits_{z}h_{x,y,z}C_z.
\end{align*}

\subsection{} Let $\A^+ = \Z[q^{1/2}]$. To go further we need the following assumptions
about the Coxeter group:
\begin{itemize}
\item[--] $(W,S)$ is {\it crystallographic}, which means that for any $s\neq t$ in $S$ the
product $st$ has order $2$, $3$, $4$, $6$ or $\infty$;
\item[--] $(W,S)$ is {\it bounded}, which means that there exists an integer $N\ge 0$ such
that $q^{N/2}f_{x,y,z}\in\A^+$ for all $x,y,z\in W$, or equivalently,
$q^{N/2}h_{x,y,z}\in\A^+$ for all $x,y,z\in W$.
\end{itemize}

\subsection{}
Right-angled Coxeter groups are crystallographic. We prove that they are bounded:

\begin{lemma}
If in right-angled Coxeter group $(W,S)$\ \
$ts_1\dots s_n = \hat ts_1\dots \hat s_j\dots s_n$ ($t, s_i\in S$,
$i = 1, \dots, n$), then $t = s_j$ and $t$ commutes with $s_1,\dots ,s_{j-1}$.
\end{lemma}

\begin{proof}
By \cite{Tits} there is a sequence of elementary M-operations taking word $ts_1\dots s_n$
to $\hat ts_1\dots \hat s_j\dots s_n$ ($t, s_i\in S$). In case of right-angled groups
the operations are
\begin{itemize}
\item[(i)] $ss \to 1$;
\item[(ii)] $st \to ts$
\end{itemize}
($s,t \in S$). Both operations preserve parity of the number of a simple reflection
occurrences in the word, which follows $t = s_j$. We shall prove that $t$ commutes
with $s_1,\dots ,s_{j-1}$ by induction on $j$. If $j-1 = 1$, then we have $ts_1t = s_1$;
$ts_1 = s_1t$. Now, if $ts_1\dots s_n$ is equivalent to
$\hat ts_1\dots s_{j-1}\hat ts_{j+1}\dots s_n$ and $j>2$, then there is an M-operation
which decreases the number of simple reflections between two $t$'s. It can be
only an operation of type~(ii), so $t$ commutes with some $s_i$ with $i\in \{1,\dots ,j-1\}$.
It follows by the induction hypothesis that $t$ also commutes with the
remaining $s_i$.
\end{proof}

\begin{theorem} Let $(W,S)$ be a right-angled Coxeter group. Then it is bounded with
$N = max\ l(w_0^{S'})$, where $S'$ runs through all subsets of $S$ such that corresponding
subgroup $(W^{S'}\!, S')$ is finite and $w_0^{S'}$ denotes the longest element of $W^{S'}$.
\end{theorem}

\begin{proof}
Take arbitrary $x, y \in W$ and let $s_k\dots s_1$, $t_1\dots t_n$ be corresponding
reduced expressions ($s_i, t_j \in S$, $i=1,\dots ,k$, $j=1,\dots ,n$). We have
$$\widetilde T_x \widetilde T_y = \widetilde T_{s_k}\dots \widetilde T_{s_1}
\widetilde T_y.$$
If $l(s_1y) = 1 + l(y)$, then $\widetilde T_{s_1} \widetilde T_y = \widetilde T_{s_1y}$;
otherwise, $l(s_1y) = l(y) - 1$ and
$\widetilde T_{s_1} \widetilde T_y = (q^{1/2} - q^{-1/2})\widetilde T_y + \widetilde T_{s_1y}$.

Proceeding inductively we obtain (similarly to \cite{L1a}):

$$\widetilde T_x \widetilde T_y = \sum_{I} (q^{1/2} - q^{-1/2})^{p_I}\widetilde T_I$$
where $I$ ranges over all subsets $i_1<\dots <i_{p_I}$ of $\{1,\dots ,k\}$ such that
$$ s_{i_l} \dots \hat s_{i_{l-1}} \dots \hat s_{i_1} \dots s_1y <
\hat s_{i_l} \dots \hat s_{i_{l-1}} \dots \hat s_{i_1} \dots s_1y $$
for $l = 1,\dots ,p_I$, and
$\widetilde T_I = s_k \dots \hat s_{i_p} \dots \hat s_{i_1} \dots s_1y$,  $p = p_I = \#I.$

From the lemma it follows that for any such $I$ there exists a reduced expression of $x$
in which $i_j = j$ and $s_{i_j} \in \L(y)$ for $j = 1,\dots ,l.$
Denoting by $\Gamma_y$ the subgroup of $W$ generated by simple reflections from $\L(y)$,
we have $s_{i_p}\dots s_{i_1}$ is a reduced expression of an element from $\Gamma_y$.
To complete the proof it remains to show that for any $y\in W$ subgroup $\Gamma_y\le W$
is finite.
Really, let $s\neq t\in \L(y)$. This means that $y$ has a reduced expression
$su_2\dots u_n$ ($s, u_j \in S$) and
$tsu_2\dots u_n = \hat tsu_2 \dots \hat u_i \dots u_n$, so again by the lemma $st = ts$,
which follows $\Gamma_y\cong \Z_2^p$ ($p = \#\L(y)$) is a finite group.
\end{proof}

\subsection{} \label{a(z)}
Let $a(z)$ be the smallest integer such that $q^{a(z)/2}h_{x,y,z}\in\A^+$ for
any $x,y\in W$. It follows that $0\le a(z)\le N$ for any $z\in W$. Here are some important
properties of the function $a$ obtained in \cite{L1a}, \cite{L2a}:

\medskip

{\noindent\normalsize {\sc Properties of $a(z)$:}}
{\it
\begin{itemize}
\item[(i)] $a(w) = a(w^{-1})$ for all $w\in W$;
\item[(ii)] $a(w) = 0$ if and only if $w = e$;
\item[(iii)] the function $a$ is constant on the two-sided cells of $W$;
\item[(iv)] $a(w) \le l(w)$ for all $w\in W$.
\end{itemize}
}

\subsection{}\label{distinguished}
Define a subset $\D\subset W$ as follows:
$$\D = \{z\in W \mid a(z) = l(z) - 2\delta(z)\},$$
where $\delta(z)$ is the degree of $P_{e,z}$ as a polynomial in $q$.

It can be shown that $d^2 = e$ for any $d\in\D$. The elements of $\D$ are called
{\it distinguished involutions} of $W$. The following theorem was proved in \cite{L2a}:

\begin{theorem}
Any left cell contains a unique $d\in\D$.
\end{theorem}

We shall use this powerful result to distinguish left cells of the
groups $P_n$ in the next section.

\section{Cells}

We mainly consider left cells and corresponding left $\H$-modules. The results concerning
right cells are entirely similar. Two-sided cells and $\H$-bimodules are obtained via
combination of the left and right-sided ones.

Let us first suppose that $(W,S)$ is an arbitrary Coxeter group. We shall always use $s,t$
(possibly with subscribed indexes which will be not connected with the initial ordering of
generators in the case of the group $P_n$) to denote the elements from $S$.

\subsection{}\label{sec51}

Let $w_1, w_2 \in W$. We say that $w_1$ and $w_2$ belong to the same {\it left precell} and
write $w_1 \sim_{lp} w_2$ if there exists $w_L \ne 1$ such that $w_1$, $w_2$ have reduced
expressions $w_1'w_L$, $w_2'w_L$, respectively, and
\begin{itemize}
\item[(a)] if $l(w_L) > 1$, then $w_1'$, $w_2'$ are either trivial or segments in $w_1$, $w_2$,
resp.;
\item[(b)] if $l(w_L) = 1$, then $w_1$ and $w_2$ are lines (with $\R(w_1) = \R(w_2) = w_L$).
\end{itemize}
Let us also suppose $1\sim_{lp} 1$.

\begin{prop}  Relation $\sim_{lp}$ has the following properties:
\begin{itemize}
\item[(i)] it is an equivalence relation on $W$;
\item[(ii)] each left precell $\Gamma$ contains a unique shortest element $w_L = w_L(\Gamma)$
such that $\Gamma = \{w \mid w = w'w_L\}$ where $w'$ is either trivial, or a segment in $w$,
or $w'w_L$ is a line (as in the definition of $\sim_{lp}$).
\end{itemize}
\end{prop}

\begin{proof} To prove (i) we need only to check the transitivity of $\sim_{lp}$. Suppose
$$w_1 \sim_{lp} w_2,\ w_2 \sim_{lp} w_3.$$
Then we have
$$w_1 = w_1'w_{L1},\ w_2 = w_2'w_{L1} = w_2''w_{L2},\ w_3 = w_3'w_{L2},$$
with $w_1'$, $w_2'$, $w_2''$, $w_3'$ as in the definition of $\sim_{lp}$. We want
to show $w_1\sim_L w_3$ and we can suppose that $w_1$, $w_2$, $w_3$ are all unequal
since otherwise it is trivial. We now need to do some routine case by case considerations:
\newline
1) $l(w_{L1}) = 1$. Then $w_2$ is a line, so $w_2''$ cannot be its proper segment; consequently
we have only two possibilities:
\begin{itemize}
\item[a)]
$l(w_{L2}) = 1$; $w_{L2} = w_{L1}$; $w_3 \sim_L w_1$.
\item[b)]
$w_2'' = 1$; $w_3 = w_3'w_{L2} = w_3'w_2$, $w_3'\neq 1$ is a segment in $w_3$ but this is
impossible by the definition of segment since $w_2$ is a line.
\end{itemize}
2) The case $l(w_{L2}) = 1$ is entirely similar.
\newline
3)  $l(w_{L1}) > 1$, $l(w_{L2}) > 1$. There are again two possibilities:
\begin{itemize}
\item[a)]
$w_2' = 1$; $w_1 = w_1'w_{L1} = w_1'w_2 = w_1'w_2''w_{L2}$, $w_1'\neq 1$ is a segment in
$w_1$. We have either $w_2'' = 1$, which follows $w_3 \sim_L w_1$, or $w_2''$ is a segment in
$w_2''w_{L2}$ which leads to a contradiction with the definition of segment.
\item[b)]
$w_2'$ is a segment in $w_2$ then either $w_2' = w_2''$ and $w_3 \sim_L w_1$ or
$w_2'' = 1$; $w_3 = w_3'w_{L2} = w_3'w_2 = w_3'w_2'w_{L1}$, which is impossible by the
definition of segment.
\end{itemize}
So in all the possible cases we obtain $w_1\sim_L w_3$.

\medskip

To prove (ii) we can take for $w_L$ (any) shortest element of the equivalence class $\Gamma$.
Then any $w$ in $\Gamma$ will have the required form by the definition and uniqueness of
$w_L$ follows. It is also possible to deduce the existence and uniqueness of $w_L$ from the
proof of (i).
\end{proof}

The existence of the canonical representatives $w_L$ of the left precells implies
that the following definitions make sense:

{\it Dimension} of a non-unit left precell $\Gamma$ $\dim(\Gamma) = \#\L(w_L)$, since
$l(w_L) \geq 1$ we have $\dim(\Gamma) \geq 1$.  Given two left precells $\Gamma_1$,
$\Gamma_2$ we say that $\Gamma_1 \le_L \Gamma_2$,  $\Gamma_1 \prec \Gamma_2$,
$\Gamma_1 - \Gamma_2$ and $\Gamma_1 \sim_L \Gamma_2$ if $w_L(\Gamma_1) \le_L w_L(\Gamma_2)$,
$w_L(\Gamma_1) \prec w_L(\Gamma_2)$,  $w_L(\Gamma_1) - w_L(\Gamma_2)$ and
$w_L(\Gamma_1) \sim_L w_L(\Gamma_2)$, respectively. We also define
$\L(\Gamma) = \L(w_L(\Gamma))$ and $\R(\Gamma) = \R(w_L(\Gamma))$.

\subsection{}
\begin{lemma} A left cell in $W$ is a union of the left precells which are $\sim_L$-equivalent
to each other.
\end{lemma}

\begin{proof}
By Lemma \ref{lemma1} if $w_1 \sim_{lp} w_2$, then $w_1 \sim_L w_2$. The corresponding
chains $x_0 - x_1 - \dots - x_k$, joining $w_1$  ($w_2$) with $w_L$ and having
property $\L(x_i)\cap \L(x_{i+1}) = \emptyset$ for any $i$ (which is actually stronger
than it is required for $\sim_L$), are obtained from the lines $w_1'$ ($w_2'$) defined
in Proposition \ref{sec51} by the rule $x_0 = w_1'w_L = t_1\dots t_kw_L$,
$x_i = t_{i+1}\dots t_kw_L$. Since both $\sim_L$ and $\sim_{lp}$ are equivalence
relations the remaining part of the statement follows easily.
\end{proof}

\medskip

The language of precells seems to be very appropriate for the description of the cells
of right-angled Coxeter groups. We shall consider in detail the case $W = P_n$. Similar
methods can be applied to the other right-angled Coxeter groups as well, we are going to
study these cases elsewhere.

\subsection{} \label{theorem}
\begin{theorem} The non-unit left cells of the group $P_n$ ($n\ge 5$) are:
\begin{itemize}
\item[(i)] $n$ cells corresponding to the $1$-dimensional left precells of $P_n$ defined by
the generators of $P_n$;
\item[(ii)] infinitely many cells which are equivalence classes of the left precells with
the canonical representatives $\Gamma(w_L)$, such that $w_L = t_1t_2w_L'$ with
$t_1t_2 = t_2t_1$ ($t_1,t_2\in S$) and $w_L'$ is a segment in $w_L$.
\end{itemize}
\end{theorem}

\begin{proof} 1) We first show that each left precell of $P_n$ belongs to an at least
one left cell defined in the statement. If $\dim(\Gamma) = 1$, then $l(w_L(\Gamma)) = 1$
(as it easily follows from the definitions of the precell and $w_L$), so $\Gamma$
is one of the $1$-dimensional precells from $(i)$. By the definition of the group
$P_n$ the dimensions of its left precells are not greater than $2$ (this also
follows from the representation of $P_n$ as a group of isometries of the hyperbolic
plane); thus it remains to consider a left precell $\Gamma$ with $\dim(\Gamma) = 2$.

Let $w_L = w_L(\Gamma)$, we have
$w_L = t_{1,1}\dots t_{1,n_1}u_1 t_{2,1}\dots t_{2,n_2}u_2 \dots  t_{k,1}\dots t_{k,n_k}u_k$
where for any admissible $i, j$ the subwords $u_i$ are either trivial or segments in $w_L$,
$n_i\ge 2$ and $t_{i,j}t_{i,j+1} = t_{i,j+1}t_{i,j}$. Define two elementary moves between
reduced words:

\begin{itemize}
\item[A:] $t_1t_2t_3x \to t_2t_3x$ for $t_1t_2 = t_2t_1$,  $t_2t_3 = t_3t_2$ and arbitrary
subword $x$;
\item[B:] $s_1s_2ut_1t_2x \to t_1t_2x$ for $t_1t_2 = t_2t_1$, $s_1s_2 = s_2s_1$, $u$ is a
segment or $u=1$ and arbitrary $x$.
\end{itemize}
By applying this elementary moves to $w_L$ one can obtain the word $t_{k,n_k-1}t_{k,n_k}u_k$
($u_k$ is a segment) which defines a canonical precell in (ii). We shall show that the moves
produce $\sim_L$-equivalent words.

The equivalence $w = t_1t_2t_3x \sim_L w_0 = t_2t_3x$ is easy: we have $w = t_1w_0$,
$\L(w_0) = \{t_1, t_2\}$, $\L(w) = \{t_2, t_3\}$ and $t_3 \neq t_1$ (because all the
expressions are reduced), so $w \succ w_0$, $\L(w_0)\not\subset \L(w)$ and $\L(w)\not\subset \L(w_0)$.

To prove that move B is a left equivalence we shall use a supplementary construction.
Having $w = s_1s_2ut_1t_2x$, $w_0 = t_1t_2x$, define $w^* = t_1u^{-1}s_1s_2ut_1t_2x$.
Note that $u$ is a segment in $w$ implies $u^{-1}$ and $u$ are segments in $w^*$.
We are going to show the following relations:

\def\lel{\raisebox{2pt}{$\scriptstyle\le_L$}}
\def\gel{\raisebox{2pt}{$\ \scriptstyle\ge_L$}}
\def\sil{\raisebox{2pt}{$\ \ \scriptstyle\sim_L\ $}}
$$
\begin{CD}
w_0 & \frac{\lel}{\phantom{xxxxxxxxxxxxxxxxx}} & w^* &
\frac{\sil}{} \:w_1^*  \frac{\sil}{}\dots \frac{\sil}{}\:w_l^*\frac{\sil}{} \:w\\
\Big|\: &  & \Big|\; & \\
w_1 & \frac{\gel}{}\:w_2\frac{\gel}{}\dots\frac{\gel}{} & \:w_n & \\
\end{CD}
$$

\medskip

\noindent where $w_{i+1}$ is obtained from $w_i$ by adding at the left the next letter
of $w^*$ and $w^*_{j+1}$ is obtained from $w^*_j$ by deleting a letter at the left.

The difficult part is to prove $w_0\le_L w^*$ since all the other chains are just
of the form $x - y$ with $|l(x)-l(y)| = 1$ and so satisfy the definitions (one can
also note that $w^*\sim_{lp}w$).

Let $u = u_1\dots u_k$ with $u_i\in S$ is a (the) reduced expression of $u$. It is enough to show
that the coefficient $\mu(w_0,w^*)$ of the $q^{(l(w^*)-l(w_0)-1)/2} = q^{k+1}$ in $P_{w_0,w^*}$
is not $0$. Let us make use of the following formula for the polynomials $P_{y,w}$ obtained
in the proof of existence of the $C_w$-basis in \cite{KL}:
\begin{align*}
P_{y,w} = q^{1-c}P_{sy,v} + q^cP_{y,v} -
\sum_{\textstyle{y\le z\prec v \atop sz < z }}
\mu(z,v)q_z^{-1/2}q_v^{1/2}q^{1/2}P_{y,z} \qquad (y\leq w),
\end{align*}
where $w = sv$ with $l(w) = 1 + l(v)$, $c=1$ if $sy < y$, $c=0$ if $sy>y$ and $P_{x,v} = 0$ unless
$x\leq v$.

We have
$$P_{w_0, w^*} = P_{t_1t_2x, t_1u^{-1}s_1s_2ut_1t_2x} = P_{t_2x,v} + qP_{t_1t_2x, v} - \varSigma_0$$
with $v = u^{-1}s_1s_2ut_1t_2x$ and $\varSigma_0 = \sum \mu(z,v)q_z^{-1/2}q_v^{1/2}q^{1/2}P_{y,z}$
where the summation is over $z$ such that $t_1t_2x\le z \prec v$, $t_1z < z$.

$\mu(t_2x, v) = 0$ because $\L(v) = u_k \not\in \L(t_2x)$ since $u$ is a segment in $w$,
and so by Lemma~\ref{lemma1} $P_{t_2x,v}$ has the maximal possible degree ($=k+1$)
if and only if $u^{-1}s_1s_2ut_1t_2x = u_kt_2x$, which is impossible. Consider the
second term:
$$qP_{t_1t_2x, v} = q^2P_{u_kt_1t_2x, u_{k-1}\dots u_1s_1s_2u_1\dots u_kt_1t_2x} + qP_1 -
q\varSigma_1,$$
defining for $i=1,\dots ,k:$
\begin{align*}
P_i & = P_{y_i, v_i},\cr
\varSigma_i & =
\sum\limits_{\textstyle{y_i \le z \prec v_i \atop u_{k-i+1}z < z}}
\mu(z,v)q_z^{-1/2}q_{v_i}^{1/2}q^{1/2}P_{y_i,z},\cr
y_i & = u_{k-i+2}\dots u_kt_1t_2x,\cr
v_i & = u_{k-i}\dots u_1s_1s_2u_1\dots u_kt_1t_2x
\end{align*}
(with the conventions that $u_{k-i+2}\dots u_k = 1$ for $i=1$ and $u_{k-i}\dots u_1 = 1$ for $i=k$).

The remarkable point is that the coefficients of $q^{k-i+2}$ in $qP_i$ and $\varSigma_{i-1}$
for $i = 1, \dots ,k$ are equal! Really, consider the sum $\varSigma_i$. We have
$\L(v_i) = u_{k-i} \not\in \L(z)$ (since $u_{k-i+1}\in \L(z)$ and it does not commute with $u_{k-i}$)
so by Lemma~\ref{lemma1}, $z\prec v_i$ implies $v_i = u_{k-i}z$, $\mu(z,v_i) = 1$,
$z = u_{k-i-1}\dots u_1s_1s_2u_1\dots u_kt_1t_2x = v_{i+1}$ and this is possible only if
$$u_{k-i-1} = u_{k-i+1}. \qquad (*)$$
In this case we have $P_{y_i,z} = P_{y_i, v_{i+1}} = P_{y_{i+1}, v_{i+1}}$ (the last equality
is the consequence of $\L(y_i) = u_{k-i+2} \neq u_{k-i+1} = \L(v_{i+1})$) and so
$\varSigma_i = qP_{i+1}$.
It remains to check that if $\mu(y_{i+1}, v_{i+1}) \neq 0$, then the equality $(*)$ holds, which
can be easily done by supposing on the contrary that $u_{k-i-1} \neq u_{k-i+1}$ and applying
Lemma \ref{lemma1}. The case of $\varSigma_0$ should be considered separately, of course, but
appears to be very similar.

So the leading terms of $qP_i$ and $\varSigma_{i-1}$ repetitively cancel each other and we finally
obtain that the coefficient of $q^{k+1}$ in $P_{w_0, w^*}$ is equal to the leading coefficient
of $q^{k+1}P_{u_1y_k, v_k}$, which is equal to $1$ because $u_1y_k < v_k$ and $l(v_k) - l(u_1y_k) = 2$.

This proves $w_0\prec w^*$ with $\mu(w_0, w^*) = 1$.

\medskip

2) It remains to show that the cells defined in the statement do not intersect.
We shall use the distinguished involutions.

For the non-unit elements $z\in P_n$ we have $a(z)\in\{ 1,2\}$. If $z$ is in a cell of type~(i),
then there exists $s\sim_L z$ and by the Properties~\ref{a(z)} $a(z) = a(s) = 1$. Now let $z$ is in a
cell of type~(ii). Using Moves~A, B from the first part of this proof and their right-side analogs
we see that $z$ belongs to the same two-sided cell as $st$ ($s,t\in S$, $(st)^2 = 1$). It is
easy to show that $a(st) = 2$ (take $x = y = st$ in the definition of the function $a$), so
again by \ref{a(z)}\ \ $a(z) = 2$.

It immediately follows from the definitions that $s_i\in\D$, $i = 1,\dots,n$, so we have the
distinguished involutions for each of the cells of type~(i). Now consider a cell of type~(ii)
with the represenative $\Gamma(w_L)$ as in the statement of the theorem. An element
$z = {w'_L}^{-1}t_1t_2w_L'\in\Gamma(w_L)$ is an involution, we shall see that $z\in\D$.
Suppose $w_L' = u_1\dots u_k$, $u_i\in S$. We have
$$P_{e,z} = P_{e,u_k\dots u_1t_1t_2u_1\dots u_k} = P_{u_k,u_k\dots u_1t_1t_2u_1\dots u_k}.$$
We see that the argument which was used to show that $w_0\le_L w*$
in part~(1) of the proof works without any changes in this case either and gives
$deg(P_{e,z}) = k$. So
$$ l(z) - 2\delta(z) = 2k+2 - 2k = 2 = a(z)$$
and $z\in\D$ by the definition.

It remains to apply Theorem~\ref{distinguished} to distinguish all the left cells.
\end{proof}

\begin{figure}[ht]
\psfig{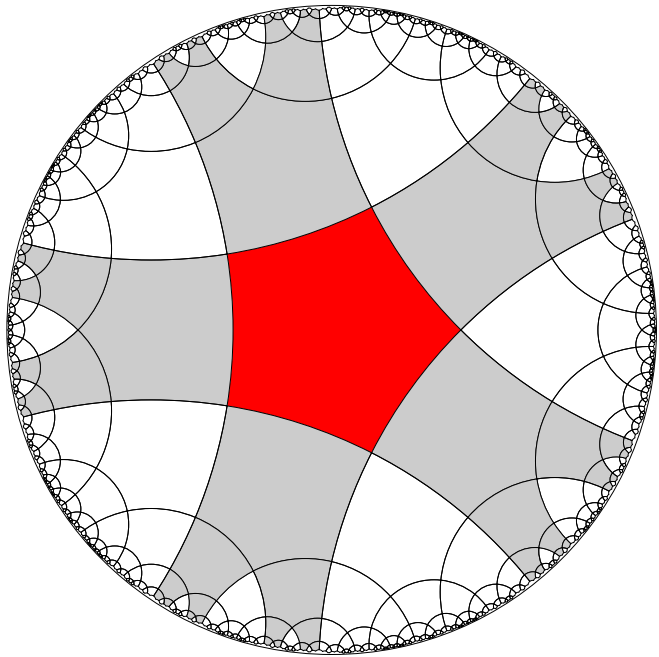}
\caption{Cells of the group $P_5.^1$}
\vskip5mm
\end{figure}

\subsection{} It was pointed out to me by V.~Ostrik that the cells of type~(i) were previously
considered in \cite{L1}. There the cells and corresponding representations were constructed for
an arbitrary Coxeter group and then thoroughly studied in the finite and affine cases.

\medskip

Left cells of the group $P_n$ can be visualized on the corresponding tessellation of the hyperbolic
plane. Figure~2 presents the cells for the group $P_5$: the pentagon in the center is the unit cell,
five shaded regions represent the cells of type~(i) all giving a one two-sided cell
(see Corollary~\ref{cor3}), each white region represents a cell of type~(ii) and altogether
they form the third two-sided cell.


\medskip

Below we give several corollaries from Theorem~\ref{theorem} and its proof.

\subsection{}
\begin{corollary}
The distinguished involutions of the group $P_n$ are
$$\mathcal{D} = \{1\} \cup \{s_1,\dots, s_n\} \cup
\{ustu^{-1}\mid (st)^2 = 1,\ u\ {\rm is\ a\ segment\ in }\ ustu^{-1}\}.$$
\end{corollary}

The distinguished involutions are related with algebra $\J$ defined in~\cite{L2a} which may
be regarded as an asymptotic version of Hecke algebra $\H$. Using the methods from~\cite{L2a}
this corollary can be applied to retrieve a partial structure information about algebra $\J$
of the group $P_n$.

\subsection{}
\footnotetext[1]{\ This picture uses W.~Casselman's PostScript library for hyperbolic geometry.}
\begin{corollary}
$W$-graphs associated to the left cells of type~(i) are infinite rooted trees (binary trees
for the group $P_5$, see Figure~3), while the $W$-graphs associated to the type~(ii) cells
admit infinitely many different cycles.
\end{corollary}

\begin{figure}[ht]
\psfig{file=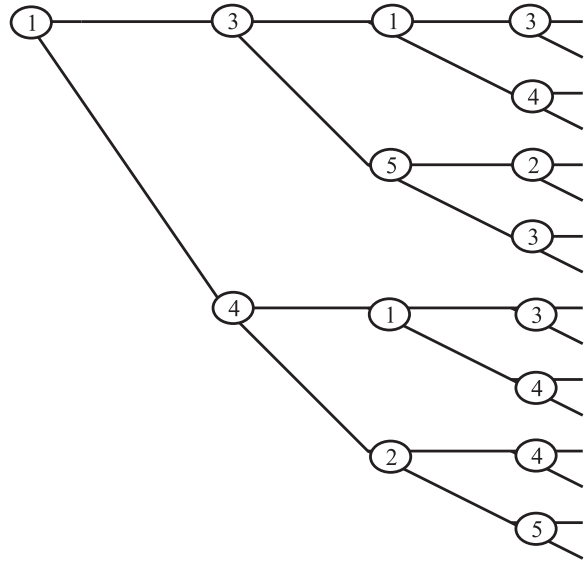, scale=.75}
\caption{Example of a $W$-graph corresponding to a left cell of type~(i) of the group $P_5$
(all $\mu(x,y) = 1$, the vertices are represented by circles with the corresponding subsets
of $S$ inside).}
\end{figure}

One can see that all the $W$-graphs corresponding to the cells of type~(i) of the group $P_n$ define
equivalent representations of $\H$, with the equivalences induced by the cyclic permutations of the
simple reflections $s_i\in S$. We suppose that the representations corresponding to the cells of
type~(ii) are also all equivalent, but this does not readily follow from the above arguments.

\subsection{}\label{cor3}
\begin{corollary} The partition of the group $P_n$ ($n\ge 5$) onto two-sided cells consits
of $3$ elements:
\begin{itemize}
\item[-] the unit cell corresponding to the trivial representation of $H$;
\item[-] the union of the left cells of type (i) (a $1$-dimensional cell), the corresponding
$W$-graph is an infinite tree;
\item[-] the union of the left cells of type (ii) (a $2$-dimensional cell), the corresponding
$W$-graph admits infinitely many different cycles.
\end{itemize}
\end{corollary}

The remarkable point about this corollary is that to establish it we actually need only
the Moves A, B from the proof \ref{theorem} with their right-side analogs, and so we do
not use part (2) of the argument which relies on certain very strong results about
distinguished involutions.


\end{document}